\documentclass{article}
\usepackage{amsmath}

\newtheorem{theorem}{Theorem}

\newtheorem{lemma}[theorem]{Lemma}

\newenvironment{proof}[1][Proof]{\noindent\textbf{#1.} }{\ \rule{0.5em}{0.5em}}
\input{tcilatex}

\begin{document}

\author{Alan Horwitz \\
Penn State University\\
25 Yearsley Mill Rd.\\
Media, PA 19063\\
alh4@psu.edu}
\date{8/11/03}
\title{The locus of centers of ellipses inscribed in quadrilaterals}
\maketitle

{\Large Introduction}

\qquad Let $R$ be a four--sided \textbf{convex} polygon in the $xy$ plane. A
problem often referred to in the literature as Newton's problem, was to
determine the locus of centers of ellipses inscribed in $R$. By inscribed we
mean that the ellipse lies inside $R$ and is tangent to each side of $R$.
Chakerian(\cite{CH}) gives a partial solution of Newton's problem using
orthogonal projection, which is the solution actually given by Newton, which
we state as

\begin{theorem}
\label{T1}Let $M_{1}$ and $M_{2}$ be the midpoints of the diagonals of $R$.
Then if $E$ is an ellipse inscribed in $R,$ the center of $E$ must lie on $Z$%
, the open line segment connecting $M_{1}$ and $M_{2}$.
\end{theorem}

However, Theorem \ref{T1} does not really give the precise locus of centers
of ellipses inscribed in $R$. It is stated in (\cite{DO}, pp. 217--219) that
the locus of centers of ellipses inscribed in $R$ actually \textbf{equals} $%
Z $, but Newton only proved that the center of $E$ must lie on $Z$, as is
noted in (\cite{CH}). Indeed, it is not even clear that an ellipse \textbf{%
exists} which is inscribed in $R,$ let alone whether \textbf{every point} of 
$Z$ is the center of such an ellipse. The main result of this note is that
it is indeed the case that \textbf{every point} of $Z$ is the center of an
ellipse inscribed in $R$. This result was actually proved by the author in (%
\cite{HR}, Theorem 11), but the approach given here is decidedly different
and much shorter and more succinct. In addition, we are also able to prove
that there is a unique ellipse of maximal area inscribed in $R$. While it is
perhaps possible to prove these results using orthogonal projection, we use,
instead, a theorem of Marden(\cite{MD}, Theorem 1) relating the foci of an
ellipse tangent to the lines thru the sides of a triangle and the zeros of a
partial fraction expansion. We state the part we shall use here.

\begin{theorem}
\label{T2}(Marden): Let $F(z)=\dfrac{t_{1}}{z-z_{1}}+\dfrac{t_{2}}{z-z_{2}}+%
\dfrac{t_{3}}{z-z_{3}}$, $t_{1}+t_{2}+t_{3}=1,$ and let $Z_{1}$ and $Z_{2}$
denote the zeros of $F(z)$. Let $L_{1},L_{2},L_{3}$ be the line segments
connecting $z_{2},z_{3}$, $z_{1},z_{3}$, and $z_{1},z_{2},$ respectively. If 
$t_{1}t_{2}t_{3}>0,$ then $Z_{1}$ and $Z_{2}$ are the foci of an ellipse, $E$%
, which is tangent to $L_{1},L_{2},$ and $L_{3}$ in the points $\zeta
_{1},\zeta _{2},\zeta _{3}$, where $\zeta _{1}=\dfrac{t_{2}z_{3}+t_{3}z_{2}}{%
t_{2}+t_{3}}$, $\zeta _{2}=\dfrac{t_{1}z_{3}+t_{3}z_{1}}{t_{1}+t_{3}}$, $%
\zeta _{3}=\dfrac{t_{1}z_{2}+t_{2}z_{1}}{t_{1}+t_{2}}$, respectively.
\end{theorem}

{\Large Main Result}

\begin{theorem}
\label{main}\qquad Let $R$ be a four--sided \textbf{convex} polygon in the $%
xy$ plane and let $M_{1}$ and $M_{2}$ be the midpoints of the diagonals of $%
R $. Let $Z$ be the open line segment connecting $M_{1}$ and $M_{2}$. If $%
(h,k)\in Z$ then there is a unique ellipse with center $(h,k)$ inscribed in $%
R$.
\end{theorem}

We shall now prove Theorem \ref{main} for the case when no two sides of $R$
are parallel. Such a quadrilateral is sometimes called a trapezium. Our
methods extend easily to the case when exactly two sides of $R$ are
parallel, that is, when $R$ is a trapezoid. Of course, if $R$ is a
parallelogram, then the midpoints of the diagonals coincide, and the line
segment $Z$ is just a point. Since ellipses, tangent lines to ellipses, and
four--sided convex polygons are preserved under affine transformations, we
may assume that the vertices of $R$ are $(0,0),(1,0),$ $(0,1)$, and $(s,t)$
for some real numbers $s$ and $t$. Let $I$ denote the open interval between $%
\dfrac{1}{2}$ and $\dfrac{1}{2}s$. Then $M_{1}=\left( \dfrac{1}{2},\dfrac{1}{%
2}\right) ,M_{2}=\left( \dfrac{1}{2}s,\dfrac{1}{2}t\right) $, and the
equation of the line thru $M_{1}$ and $M_{2}$ is 
\begin{equation*}
y=L(x)=\dfrac{1}{2}\dfrac{s-t+2x(t-1)}{s-1},x\in I
\end{equation*}%
Since $R$ is convex, four--sided and no two sides of $R$ are parallel, it
follows easily that 
\begin{equation*}
\text{ }s>0,t>0,s+t>1,\text{ and }s\neq 1\neq t
\end{equation*}

We shall need the following lemmas.

\begin{lemma}
\label{L2}If $h\in I$ and $s+t>1$, then $s+2h(t-1)>0$
\end{lemma}

\begin{proof}
If $t>1$, then $s,h,$ and $t-1$ are all positive. If $t\leq 1$ and $s\geq 1$%
, then $I=\left( \dfrac{1}{2},\dfrac{1}{2}s\right) \Rightarrow s+2h(t-1)\geq
s-2h>0$. Finally, if $t\leq 1$ and $s\leq 1$, then $I=\left( \dfrac{1}{2}s,%
\dfrac{1}{2}\right) \Rightarrow s+2h(t-1)>s+t-1>0$.
\end{proof}

We leave the proof of the next lemma to the reader.

\begin{lemma}
\label{L4}Let $E_{1}$ and $E_{2}$ be ellipses with the same foci. Suppose
also that $E_{1}$ and $E_{2}$ pass through a common point, $z_{0}$. Then $%
E_{1}=E_{2}$.
\end{lemma}

\textbf{Proof of Theorem \ref{main}: }Let $L_{1}$: $y=0,L_{2}$: $x=0,L_{3}$: 
$y=\dfrac{t}{s-1}(x-1)$, and $L_{4}$: $y=1+\dfrac{t-1}{s}x$ denote the lines
which make up the boundary of $R$. $L_{1},L_{2},$ and $L_{3}$ form a
triangle, $T_{1}$, whose vertices are the complex points $z_{1}=0$, $z_{2}=1$%
, and $z_{3}=-\dfrac{t}{s-1}i$. $L_{1},L_{2},$ and $L_{4}$ form a triangle, $%
T_{2}$, whose vertices are the complex points $w_{1}=0$, $w_{2}=i$, and $%
w_{3}=-\dfrac{s}{t-1}$. First, we want to find ellipses $E_{1}$ and $E_{2}$
tangent to $L_{1},L_{2},$ and $L_{3},$ and to $L_{1},L_{2},$ and $L_{4}$,
respectively. We shall use Theorem \ref{T2}, so that $E_{1}$ has foci $Z_{1}$
and $Z_{2}$, which are the zeros of $F(z)=\dfrac{t_{1}}{z}+\dfrac{t_{2}}{z-1}%
+\dfrac{t_{3}}{z+\frac{t}{s-1}i},$ and $E_{2}$ has foci $W_{1}$ and $W_{2}$,
which are the zeros of $G(z)=\dfrac{s_{1}}{z}+\dfrac{s_{2}}{z-i}+\dfrac{%
1-s_{1}-s_{2}}{z+\frac{s}{t-1}}$. To guarantee that $E_{1}$ and $E_{2}$ are
ellipses, we require, by Theorem \ref{T2}, that $s_{1}s_{2}s_{3}>0$ and $%
t_{1}t_{2}t_{3}>0$, where $s_{3}=1-s_{1}-s_{2}$ and $t_{3}=1-t_{1}-t_{2}$.
For example, let $s=3$, $t=2$, $t_{1}=-\dfrac{1}{4}$, $t_{2}=\dfrac{3}{2}$, $%
s_{1}=\dfrac{1}{3}$, and $s_{2}=\dfrac{1}{2}$. Then $t_{1}t_{2}t_{3}=%
\allowbreak \dfrac{3}{32}>0$ and $s_{1}s_{2}s_{3}=\allowbreak \dfrac{1}{36}%
>0 $. The foci of $E_{1}$ are approximately $\allowbreak
Z_{1}=-.1957-.04\allowbreak 96i$ and $Z_{2}=-.304\,3-1.\,\allowbreak 2004i$.
Note that $E_{1}$ is \textbf{not} \textbf{inscribed} in $T_{1}$since not all
of the $t_{j}$'s are positive(see Figure 1). $\allowbreak $The foci of $%
E_{2} $ are approximately $W_{1}=-.01\allowbreak 59+.401\,9i$ and $%
W_{2}=-2.\,\allowbreak 484\,1+.09\allowbreak 81i$. Note that $E_{2}$ \textbf{%
is} \textbf{inscribed} in $T_{2}$ since all of the $s_{j}$'s are
positive(see Figure 2). Assume now that $(h,k)\in Z$, or equivalently, that $%
k=L(h),h\in I$. We want $E_{1}$ and $E_{2}$ each to have center $(h,k)$. The
center, $C_{1}$, of $E_{1}$ is $\dfrac{1}{2}\left( Z_{1}+Z_{2}\right) $. A
simple computation shows that $C_{1}=-\dfrac{1}{2(s-1)}\left(
it(t_{1}+t_{2})+(s-1)(t_{2}-1)\right) $, which, upon taking real and
imaginary parts yields $C_{1}=\left( \allowbreak \dfrac{1}{2}-\dfrac{1}{2}%
t_{2},-\dfrac{1}{2}t\dfrac{t_{1}+t_{2}}{s-1}\right) $. Similarly, the center
of $E_{2}$ is $C_{2}=\left( -\dfrac{1}{2}s\dfrac{s_{1}+s_{2}}{t-1},-\dfrac{1%
}{2}\left( s_{2}-1\right) \right) $. We actually do not require these
explicit formulas for $C_{1}$ and $C_{2}.$However, solving $(h,k)=\left(
\allowbreak \dfrac{1}{2}-\dfrac{1}{2}t_{2},-\dfrac{1}{2}t\dfrac{t_{1}+t_{2}}{%
s-1}\right) $ for $t_{1}$ and $t_{2}$ shows that the center of $E_{1}$ is $%
(h,k)$ if and only if

\begin{equation}
t_{1}=2h-1-2k\dfrac{s-1}{t}\text{, }t_{2}=1-2h  \label{e1}
\end{equation}

Similarly, solving $(h,k)=\left( -\dfrac{1}{2}s\dfrac{s_{1}+s_{2}}{t-1},-%
\dfrac{1}{2}\left( s_{2}-1\right) \right) $ for $s_{1}$ and $s_{2}$ shows
that the center of $E_{2}$ is $(h,k)$ if and only if 
\begin{equation}
s_{1}=2k-1-2h\dfrac{t-1}{s}\text{, }s_{2}=1-2k  \label{e2}
\end{equation}

So given $(h,k)\in Z$, let $s_{1},s_{2},t_{1},t_{2}$ be \textit{defined} by (%
\ref{e1}) and (\ref{e2}). Substituting $k=L(h)$ into (\ref{e1}) and (\ref{e2}%
) yields $t_{1}t_{2}t_{3}=\left( s+2h(t-1)\right) \dfrac{\left( s-2h\right)
^{2}\left( 2h-1\right) ^{2}}{t^{3}}>0$ since $h\in I$ and by Lemma \ref{L2}.
Similarly, $s_{1}s_{2}s_{3}=\left( s+2h(t-1)\right) \left( 2h-1\right)
\left( s-2h\right) \dfrac{\left( t-1\right) ^{2}}{s^{2}\left( s-1\right) ^{2}%
}\allowbreak >0$, again since $h\in I$ and by Lemma \ref{L2}. Thus,
corresponding to each $(h,k)\in Z$, we have found ellipses $E_{1}$ and $%
E_{2} $, with $E_{1}$ tangent to $L_{1},L_{2},$ and $L_{3}$, and $E_{2}$
tangent to $L_{1},L_{2},$ and $L_{4}$. However, we require \textbf{one}
ellipse, with center $(h,k)$, which is tangent to \textbf{all four lines} $%
L_{1},L_{2},L_{3}$, and $L_{4}$. \ Well, the foci of $E_{1}$ are the zeroes
of the numerator of $F(z),$ which is the polynomial 
\begin{gather*}
p(z)=\allowbreak \left( s-1\right) z^{2}+\left(
it(t_{1}+t_{2})+(s-1)(t_{2}-1)\right) \allowbreak z-it_{1}t \\
=\left( s-1\right) (z-Z_{1})(z-Z_{2})
\end{gather*}%
Similarly, the foci of $E_{2}$ are the zeros of the numerator of $G(z)$,
which is the polynomial 
\begin{gather*}
q(z)=\allowbreak \left( t-1\right) z^{2}+\left( s(s_{1}+s_{2})+\allowbreak
i\left( s_{2}-1\right) \left( t-1\right) \right) \allowbreak z-is_{1}s \\
=\left( t-1\right) (z-W_{1})(z-W_{2})
\end{gather*}

Using $k=L(h),$ (\ref{e1}), and (\ref{e2}), \ $\dfrac{p(z)}{s-1}=\dfrac{q(z)%
}{t-1}=z^{2}-2(h+iL(h))z+i\dfrac{s-2h}{s-1}$. Since $\dfrac{p(z)}{s-1}$ and $%
\dfrac{q(z)}{t-1}$ have the \textbf{same} coefficients, $E_{1}$ and $E_{2}$
have the \textbf{same}\textit{\ }\textbf{foci}. Also, by Theorem \ref{T2}, $%
E_{1}$ and $E_{2}$ are both tangent to $L_{2}$ at the point $\left( 0,\dfrac{%
1}{2}\dfrac{s-2h}{\left( s-1\right) h}\right) $. By Lemma \ref{L4}, $E_{1}$
and $E_{2}$ are identical. Hence $E=E_{1}=E_{2}$ is an ellipse, with center $%
(h,k)$, which is tangent to \textbf{all four lines} $L_{1},L_{2},L_{3}$, and 
$L_{4}$. Of course $E$ is \textbf{inscribed} in $R$ since $(h,k)\in Z\subset
R$. To prove \textit{uniqueness}, if $E_{1}$ and $E_{2}$ are distinct
concentric ellipses, then, as noted in (\cite{CH}), their four common
tangents would have to form a parallelogram. If $R$ is not a parallelogram,
then this is a contradiction. We leave the proof when exactly two sides of $%
R $ are parallel to the reader.

{\Large Maximal Area}

\qquad We now want to minimize and/or maximize the area of an ellipse
inscribed in a four--sided \textbf{convex} polygon, $R$. First we require a
generalization of a result which appears in (\cite{CH}) on the area of an
ellipse inscribed in a triangle. Chakerian's result assumes that the point $P
$ lies \textbf{inside} $ABC$, the triangle with vertices $A,B,$ and $C$,
while our result assumes that $P$ lies \textbf{outside} $ABC$. In that case,
area$\left( ABC\right) =$ area$\left( CPA\right) +$ area$\left( APB\right) -$%
area$\left( BPC\right) $. The details of the proof are similar.

\begin{lemma}
\label{L3}Given a triangle $ABC$ and a point $P\notin \partial \left(
ABC\right) $, let $\alpha =$ area$(BPC),\beta =$ area$(CPA),$ and $\gamma =$
area$(APB)$. Let $L_{1},L_{2},$ and $L_{3}$ be the three lines thru the
sides of $ABC$, and let $E$ be an ellipse with center $P$ which is tangent
to $L_{1},L_{2},$ and $L_{3}$. If $\sigma =\dfrac{1}{2}\left( \alpha +\beta
+\gamma \right) $, then area$(E)=\dfrac{4\pi }{\text{area}(ABC)}\sqrt{\sigma
\left( \sigma -\alpha \right) \left( \sigma -\beta \right) \left( \sigma
-\gamma \right) }$
\end{lemma}

Now let $A_{E}=$ area of an ellipse $E$ inscribed in $R$. We want to
maximize and/or minimize $A_{E}$ as a function of $h$, where $(h,L(h))$
denotes the center of $E$. We discuss the case when no two sides of $R$ are
parallel. Let $A=(0,0),B=(1,0),C=\left( 0,-\dfrac{t}{s-1}\right) $ , which
are the vertices of the triangle we earlier called $T_{1}$. Then area$(ABC)=%
\dfrac{1}{2}\dfrac{t}{\left\vert s-1\right\vert }$, and since $E$ is
inscribed in $ABC$, we can apply Lemma \ref{L3}, with $P=(h,k)$.
Substituting $k=L(h)$ yields $\sigma (\sigma -\alpha )(\sigma -\beta
)(\sigma -\gamma )\allowbreak =\dfrac{1}{256}t^{2}\left( -1+2h\right) \left(
s+2ht-2h\right) \dfrac{s-2h}{\left( s-1\right) ^{4}}$. By Lemma \ref{L3}, $%
A_{E}=\dfrac{\pi }{2\left\vert s-1\right\vert }\sqrt{\left( 2h-1\right)
\left( s+2h(t-1)\right) \left( s-2h\right) }$. Thus we want to optimize $%
A(h)=\left( s-2h\right) \left( 2h-1\right) (s+2h(t-1)),h\in I$. Now $%
A(1/2)=A(s/2)=0,$ and $A(h)\geq 0$ for $h\in I$ by Lemma \ref{L2}. Hence $%
A^{\prime }(h_{0})=0$ for some $h_{0}\in I$ with $A(h_{0})$ a local maximum,
and $A(h)$ dooes not attain its global minimum on $I$. Also, $A(h_{0})$ must
be the \textbf{only} local maximum of $A(h)$ on $I$, else $A^{\prime }(h)$
would have \textbf{three} zeros in $I$. Thus $A(h_{0})$ is the global
maximum of $A(h)$ on $I$. Since ratios of areas of ellipses are preserved
under affine transformations, we have proven

\begin{theorem}
Let $R$ be any given four--sided \textbf{convex} polygon in the $xy$
plane.Then there is a unique ellipse of maximal area inscribed in $R$. There
is no ellipse of minimal area inscribed in $R$.
\end{theorem}

\textbf{Example:} Take $s=4$, $t=2$, so that $R$ has vertices $(0,0),(1,0),$ 
$(0,1)$, and $(4,2)$. Then the maximal area ellipse has center $\left( 
\dfrac{4}{3},\dfrac{7}{9}\right) $.

{\Large Hyperbolas}

Using our earlier notation, let $X$ be the open line segment which is the
part of $L$ lying inside $R$, where $L$ is the line thru the midpoints of
the diagonals. If $(h,k)\in X-Z-M_{1}-M_{2}$, it is natural to think that
there should be a \textit{hyperbola}, $H$, with center $(h,k)$, which is
tangent to each line making up the boundary of $R$. This is actually
correct, but only if one considers an asymptote of $H$ to be tangent to $H$%
(at infinity, of course). \footnote{%
This was also proven in (\cite{HR}), but the statement there is not quite
correct since this author omitted the case where the"tangent line" is an
asymptote.}This is not hard to prove using the methods of this paper. An
asymptote of $H$ can arise when employing Theorem \ref{T2} since it is
possible for one of $t_{i}+t_{j},j\neq i,$ to be $0$.

\end{document}